\documentclass[11pt]{article}

\usepackage{amsfonts,amssymb,amsbsy, amsmath, dsfont}
\usepackage[german,english]{babel}
\usepackage[utf8]{inputenc}
\usepackage{graphicx}
\usepackage{epsfig}
\usepackage{gauss}
\usepackage{color}
\usepackage{amsthm}
\usepackage{framed}
\usepackage{enumerate}
\usepackage{tikz}
\usepackage{mathtools}

\usepackage[hyperindex,backref]{hyperref} %always load me last please !!!

\newcommand{\E}{{\mathbb E}}

\newcommand{\N}{{\mathbb N}}

\renewcommand{\P}{{\mathbb P}}
\newcommand{\Q}{{\mathbb Q}}
\newcommand{\R}{{\mathbb R}}

\newcommand{\Z}{{\mathbb Z}}

\newcommand{\cA}{{\cal A}}
\newcommand{\cB}{{\cal B}}
\newcommand{\cC}{{\cal C}}

\newcommand{\cE}{{\cal E}}

\newcommand{\cG}{{\cal G}}
\newcommand{\cH}{{\cal H}}
\newcommand{\cI}{{\cal I}}

\newcommand{\cN}{{\cal N}}

\newcommand{\cS}{{\cal S}}

\newcommand{\cZ}{{\cal Z}}

\newcommand{\ind}{\mathds{1}}

\DeclareMathOperator{\Vol}{Vol}

\numberwithin{equation}{section}	

\newtheoremstyle{thm}% hnamei
	{8pt}% hSpace abovei
	{8pt}% hSpace belowi
	{\itshape}% hBody fonti
	{}% hIndent amounti1
	{\bfseries }% hTheorem head fonti
	{}% hPunctuation after theorem headi
	{\newline}% hSpace after theorem headi2
	{}%

\newtheoremstyle{namedthm}% hnamei
	{8pt}% hSpace abovei
	{8pt}% hSpace belowi
	{\itshape}% hBody fonti
	{}% hIndent amounti1
	{\bfseries }% hTheorem head fonti
	{}% hPunctuation after theorem headi
	{5pt}% hSpace after theorem headi2
	{}%
	
\newtheoremstyle{def}% hnamei
	{8pt}% hSpace abovei
	{8pt}% hSpace belowi
	{}% hBody fonti
	{}% hIndent amounti1
	{\bfseries }% hTheorem head fonti
	{}% hPunctuation after theorem headi
	{\newline}% hSpace after theorem headi2
	{}%	

\newtheoremstyle{nameddef}% hnamei
	{8pt}% hSpace abovei
	{8pt}% hSpace belowi
	{}% hBody fonti
	{}% hIndent amounti1
	{\bfseries }% hTheorem head fonti
	{}% hPunctuation after theorem headi
	{5pt}% hSpace after theorem headi2
	{}%	
	
\swapnumbers
\theoremstyle{thm}
\newtheorem{Thm}{Theorem}[section]
\newtheorem{Lemma}[Thm]{Lemma}%[section]

\theoremstyle{namedthm}
\newtheorem{NamedThm}{Theorem}[section]
\theoremstyle{def}
\newtheorem{Def}[Thm]{Definition}

\theoremstyle{nameddef}

\newlength{\XWidth}

\makeatletter
\def\moverlay{\mathpalette\mov@rlay}
\def\mov@rlay#1#2{\leavevmode\vtop{%
   \baselineskip\z@skip \lineskiplimit-\maxdimen
   \ialign{\hfil$\m@th#1##$\hfil\cr#2\crcr}}}
\newcommand{\charfusion}[3][\mathord]{
    #1{\ifx#1\mathop\vphantom{#2}\fi
        \mathpalette\mov@rlay{#2\cr#3}
      }
    \ifx#1\mathop\expandafter\displaylimits\fi}
\makeatother

\begin{document}
\author{Sebastian Ziesche\thanks{Karlsruhe Institute of Technology, sebastian.ziesche@kit.edu}
}

\title{First passage percolation in Euclidean space and on random tessellations}
\date{\today}
\maketitle

\begin{abstract}
  \noindent There are various models of first passage percolation (FPP) in $\R^d$. We want to start a very general study of this topic. To this end we generalize the first passage percolation model on the lattice $\Z^d$ to $\R^d$ and adapt the results of \cite{boivin1990first} to prove a shape theorem for ergodic random pseudometrics on $\R^d$. A natural application of this result will be the study of FPP on random tessellations where a fluid starts in the zero cell and takes a random time to pass through the boundary of a cell into a neighbouring cell.

  We find that a tame random tessellation, as introduced in the companion paper \cite{ziesche2016bernoulli}, has a positive time constant. This is used to derive a spatial ergodic theorem for the graph induced by the tessellation.

  Finally we take a look at the Poisson hyperplane tessellation, give an explicit formula to calculate it's FPP limit shape and bound the speed of convergence in the corresponding shape theorem.
\end{abstract}

\maketitle

\begin{flushleft}
\end{flushleft}

\section{Introduction}\label{secintro}
First Passage Percolation was introduced by Hammersley and Welsh in 1965 as a model for the flow of a fluid through a random medium. It is introduced as a model on the graph $\cZ^d$, that is the vertex set $\Z^d$ with connections between nearest neighbours, where each edge has an i.i.d.\ passage time. The object of interest is the pseudometric $\tau$ on $\cZ^d$, which equals the time needed to travel along the shortest path (with respect to the sum of passage times of the traversed edges) from one node of the graph to another.

This concept could be applied to $\R^d$ too. Think of a tessellation $m$ where each face of a cell comes with a passage time. Then define a pseudometric $\tau(x,y)$ that minimizes the travel time over all continuous curves from $x$ to $y$, where the travel time of a curve $\gamma$ is just the sum of passage times of cell faces that where intersected by $\gamma$. This coincides with the original model, if $m$ is the cubic lattice and the passage times are i.i.d.\ .

When we first started our study as a part of the authors PhD thesis we wanted to generalize the model by replacing the fixed tessellation $m$ by some random tessellation as it was done in \cite{vahidi1990first} and \cite{vahidi1992shape} in the special case of the Poisson Voronoi tessellation.

The first aim was to prove a shape theorem. It became clear, that we would have to use or adapt the shape theorem for ergodic passage times by Boivin \cite{boivin1990first}. But when we had written down the proof, we wondered, that no particular properties of random tessellations had been used at all and we recognized that we had in fact adapted Boivins result to \emph{arbitrary random ergodic pseudometrics} on $\R^d$. The applications of this result go of course far beyond random tessellation. Models where this theorem could be applied are for example the paper on random Riemannian FPP \cite{lagatta2009shape}, the special degenerate case of Liouville FPP in \cite{ding2016upperbounds}, the FPP Model induced by the Boolean model from \cite{gouere2016positivity} or Howard and Newman's work \cite{howard1997euclidean}.

The problem with Boivin's result and our adaption is, that it doesn't rule out a degenrate limit in which case there is no real limit shape and the time constant may be zero. That means it could be possible to traverse the space at arbitrary high speed. With the application to random tessellations in mind, we found that the notion of a \emph{tame tessellation} introduced in \cite{ziesche2016bernoulli}, was perfectly fit to exclude this pathological case.

The main purpose of the present work is to lead to a better understanding of the global structure of random tessellations. In \cite{ziesche2016bernoulli} we studied the Bernoulli face percolation model on a random tessellation and we continue this work by studying FPP on random tessellations. The result we obtained is, that under reasonable regularity assumptions, a shape theorem with a non-degenerate limit shape holds for the pseudometric induced by the random tessellation with i.i.d.\ passage times.

A natural application of the shape theorem is the generalization of the multivariate ergodic theorem of Wiener \cite[Thm 10.14]{kallenberg2002foundations} to the graph induced by the random tessellation, i.e.\ the graph that has the cells of the tessellation as vertices with two vertices being adjacent if and only if they share a common $d-1$-dimensional face. This will also shed a little light on the various interpretations of Palm probabilities.

Finally we studied the model on the Poisson hyperplane tessellation as this has very nice FPP properties despite not being a tame tessellation. In fact, we are able to \emph{compute the limit shape explicitly} and bound the speed of convergence in the shape theorem.

The rest of the paper is structured in the following way. In Section \ref{sec:shape-thm-ergodic-psmetrics} we proof a weak shape theorem for random ergodic pseudometrics on $\R^d$. In Section \ref{sec:fpp-on-random-tess} we restrict ourselves to FPP on random tessellations and proof that tame tessellations satisfying certain moment conditions admit a shape theorem with a non-degenerate limit shape. Section \ref{sec:spatial-erg-thm-on-graphs} contains a very short introduction to Palm probabilities and the proof of a spatial ergodic theorem where the average is taken over random balls. Finally in Section \ref{sec:poisson-hyperplane-tess} we show, how the limit shape looks for the Poisson hyperplane tessellation and we give bounds on the deviation from the limit shape in this model.

\section{A Shape Theorem for ergodic pseudometrics}\label{sec:shape-thm-ergodic-psmetrics}
We work on an abstract probability space $(\Omega, \mathfrak{F}, \P)$ and will use various kinds of balls with respect to different metrics or pseudometrics. We write $B_r(x)$ for the Euclidean ball around $x \in \R^d$ with radius $r \geq 0$. The ball with respect to the (possibly random) pseudometric $\tau$ is denoted by $B_r^\tau(x) := \{y \in \R^d \mid \tau(x,y) \leq r\}$. On a graph $\cG = (V, E)$ with vertices $V$, $E$ and induced graph metric $d_{\cG}$ we define $B_n^{\cG}(v) := \{w \in V \mid d_\cG(v,w) \leq n\}$ for $n \in \N_0 := \N \cup \{0\}$ and $v \in V$.

A family $\tau := \{\tau(x,y) \mid x, y \in \R^d\}$ of $[0,\infty )$-valued random variables is called a \emph{random pseudometric} on $\R^d$ if a.s.
\begin{align*}
  \tau(x,x) &= 0,\\
  \tau(x,y) &= \tau(y,x),\\
  \tau(x,y) &\leq \tau(x,z) + \tau(z,y)
\end{align*}
for all $x,y,z \in \R^d$. Formally this is a random field on $\R^d \times \R^d$, i.e.\ a random element in the space $[0,\infty )^{\R^d \times \R^d}$ equipped with the usual $\sigma$-algebra that is induced by the finite dimensional projections.

The shift operator $\theta_z$, $z \in \R^d$ is defined on this space via
\[
  \theta_z (\tau(x,y)) := \tau(x-z, y-z)
\]
and a random pseudometric $\tau$ is called \emph{stationary} if $\theta_z \tau \overset{d}{=} \tau$. We denote by $\cI_x$ the $\sigma$-algebra of $\theta_x$ invariant events, i.e.\ measurable subsets $A \subset [0,\infty )^{\R^d \times \R^d}$ with $\theta_x A = A$. A random pseudometric $\tau$ is called \emph{ergodic} if $\P[\tau \in A] \in \{0,1\}$ for all $x \in \R^d$ and $A \in \cI_x$. We have to remark that this definition is a bit stronger than usual. The differences to the usual definition, where $A$ is in $\bigcap_{x \in \R^d} I_x$ are studied in \cite{pugh1971ergodic}.

We recall that a semi-norm $\mu$ is a map from $\R^d$ to $[0,\infty )$ such that
\begin{enumerate}
  \item $\mu(\lambda x) = | \lambda| \mu(x), \quad \lambda \in \R, x \in \R^d$,
  \item $\mu(x + y) \leq \mu(x) + \mu(y), \quad x,y \in \R^d$.
\end{enumerate}
This implies
\begin{align}\label{eq:lipschitzbound_on_seminorm}
  | \mu(x) - \mu(y)| \leq \max_{i \in [d]} \mu(\mathbf{e}_i) \|x - y\|_1,
\end{align}
where $\mathbf{e}_i$ with $i \in [d] := \{1, \dots, d\}$ are the standard unit vectors, and hence the Lipschitz continuity of $\mu$. In this setup, the following shape theorem, that is the continuous space pendent of the main theorem in \cite{boivin1990first}, holds. We will call this type of theorem a \emph{weak shape theorem} later on.

\begin{Thm}\label{thm:shape_theorem_random_metric}
  Let $\tau$ be a random ergodic pseudometric. If
  \begin{align}\label{eq:shape_theorem_random_metric:mom_cond}
    \E\left[ \max \left\{ \tau(0,x) \mid \|x\|_\infty \leq 1 \right\}^{d+\varepsilon} \right] < \infty
  \end{align}
  then there is a semi-norm $\mu := \mu_\tau$, such that
  \begin{align}\label{eq:shape_theorem_random_metric:uniform-limit}
    \lim_{\|x\|_2\to \infty} \frac{1}{\|x\|_2} (\tau(0, x) - \mu(x)) = 0.
  \end{align}
\end{Thm}
\emph{Proof:} Most parts of the proof look exactly like in \cite{boivin1990first} but we have to circumvent some additional problems.

For a fixed $x \in \R^d$ the limit
\begin{align}\label{eq:shape_theorem_random_metric:1}
  \mu(x) := \lim_{n\to \infty } \frac{\tau(0, nx)}{n} = \inf_{n\in\N} \frac{\E[\tau(0, nx)]}{n} = \lim_{n \to \infty} \frac{\E[\tau(0, nx)]}{n}.
\end{align}
exists due to Kingman's subadditive ergodic theorem (see \cite[Theorem 10.22]{kallenberg2002foundations} or \cite{liggett1985improved}). The triangle inequality for $\tau$ and its stationarity imply the triangle inequality for $\mu$. An easy calculation shows the homogenity of $\mu$ for rational $\lambda$ and the moment condition is more than enough to extend this to $\lambda \in \R$. Hence we only have to ensure the uniformity of the convergence in \eqref{eq:shape_theorem_random_metric:uniform-limit}.

An easy calculation shows, that \eqref{eq:shape_theorem_random_metric:mom_cond} and the Borel-Cantelli lemma imply that for any $\varepsilon > 0$ a.s.\ only a finite number of the events
\[
  \{\max \left\{ \tau(v,v+x) \mid \|x\|_\infty \leq 1 \right\} \geq \varepsilon \|v\|_2\}, \quad v \in \Z^d
\]
hold. Hence by writing $x = v + r$ with $v \in \Z^d$ and $r \in [0,1]^d$ it follows with with \eqref{eq:lipschitzbound_on_seminorm}, that it is enough to show \eqref{eq:shape_theorem_random_metric:uniform-limit} for $x \in \Z^d$.

The idea is, to show this convergence in spherical sectors
\[
  S(u, r, \delta) := \{x \in \R^d \mid \|x\|_2 \leq r,\ \langle \tfrac{x}{\|x\|_2}, u\rangle \geq 1-\delta\}
\]
with direction $u \in \cS^{d-1} := \{x \in \R^d \mid \|x\|_2 = 1\}$, length $r \in [0,\infty ]$ and opening angle $\arccos(1 - \delta) \in [0,2 \pi]$.

We define the discrete family $T := (T_e)_{e \in E(\cZ^d)}$ of passage times via
\[
  T_{\{v,w\}} := \tau(v,w), \quad \{v,w\} \in E(\cZ^d)
\]
where $E(\cZ^d) := \{\{v,w\} \subset \Z^d \mid \|v - w\|_2 = 1\}$ is the edge set of the graph $\cZ^d$. Corresponding to $T$ we define the random pseudometric $\tau'$ on $\cZ^d$ via
\[
  \tau'(v,w) := \inf_{\gamma:v \leftrightarrow w} \sum_{e \in \gamma} T_e, \quad v,w \in \Z^d
\]
where the infimum is taken over all paths in $\cZ^d$ starting in $v$ and ending in $w$. The ergodicity of $\tau$ implies that $T$ is ergodic with respect to $(\theta_v)_{v \in \Z^d}$. As we also know that
\begin{align}
  \E[T_{\{0,\mathbf{e}_i\}}^{d+\varepsilon}] \leq \E[\max\{\tau(0,x) \mid \|x\|_\infty \leq 1\}^{d+\varepsilon}] < \infty , \quad i \in [d]
\end{align}
we may apply the Maximal Lemma of \cite{boivin1990first} to obtain that there are constants $c_1, c_2 \in \R$ depending only on the dimension such that
\begin{align}
  \P\left[ \sup_{v \in \Z^d \setminus \{0\}} \frac{\tau'(0, v)}{\|v\|_2} > \lambda \right] \leq \frac{c_1 \max_{i \in [d]}\E\big[T_{\{0,\mathbf{e}_i\}}^{d+\varepsilon}\big] }{\lambda^d} = c_2 \lambda^{-d}, \quad \lambda > 0.
\end{align}
We fix a direction $u \in \cS^{d-1}$ and want to use the Maximal Lemma to show that spherical shells of any thickness contain points that are well behaved. That means, that for $\rho > 0$, $\lambda = \sqrt[d]{c_2}$ and $\delta > 0$ a.s.\ there is a large enough $r \in \R$ such that there is a $v \in (S(u,r(1+\rho),\delta) \setminus S(u,r, \delta)) \cap \Z^d$ where
\begin{align}\label{eq:shape_theorem_random_metric:3}
  A_v := \left\{ \sup_{w \in \Z^d \setminus \{v\}} \frac{\tau'(v,w)}{\|v - w\|_2} \leq \lambda \right\}
\end{align}
holds.

To this end we fix $\rho > 0$ as well as $\delta > 0$ and apply the spatial ergodic theorem (see \cite[Prop. 4.23, Example 1]{nguyen1979ergodic}) to the sequence of spherical segments $\{S(u, r, \delta) \mid r \in \N \}$ and the function $\ind_{A_0}$. We obtain that for any $\varepsilon \in (0, \P[A_0])$ and all large enough $r$
\[
  \P[A_{0}] - \varepsilon \leq \frac{1}{|S(u, r, \delta) \cap \Z^d|} \sum_{v \in S(u, r, \delta) \cap \Z^d} \ind_{A_v} \leq \P[A_{0}] + \varepsilon.
\]
Let us assume, that for a given realization of $\tau$ there is no $v$ with the property we stated in \eqref{eq:shape_theorem_random_metric:3} and above. Then
\[
  \P[A_{0}] - \varepsilon \leq \frac{1}{|S(u, r(1 + \rho), \delta) \cap \Z^d|} \sum_{v \in S(u, r, \delta) \cap \Z^d} \ind_{A_v}
\]
for large enough $r$. As $\lim_{r\to \infty } r^d/|S(u, r, \delta) \cap \Z^d| = c_3$ for some $c_3 > 0$ depending on $\delta$, we have
\begin{align*}
  (c_3 - \varepsilon) (r^d(1+\rho)^d - r^d) & \leq |(S(u,r(1+\rho),\delta) \setminus S(u,r, \delta)) \cap \Z^d | \\
  & \leq \frac{1}{\P[A_{0}] - \varepsilon} \sum_{v \in S(u, r, \delta) \cap \Z^d} \ind_{A_v} - |S(u, r, \delta) \cap \Z^d|  \\
  & \leq  |S(u, r, \delta) \cap \Z^d| \left( \frac{\P[A_{0}] + \varepsilon}{\P[A_{0}] - \varepsilon} - 1\right) \\
  & \leq (c_3 + \varepsilon) r^d \left( \frac{\P[A_{0}] + \varepsilon}{\P[A_{0}] - \varepsilon} - 1\right)
\end{align*}
for large $r$. This leads to a contradiction if $\varepsilon$ is small enough.

The existence of these well-behaved points will ensure that all points in a spherical shell are relatively close to each other. For large $r$ let $w(r) \in \Z^d$ be a point in $S(u,r(1+\rho),\delta) \setminus S(u,r, \delta)$ such that $A_{w(r)}$ holds. For $v \in S(u, \infty , \delta) \cap \Z^d$ with $\|v\|_2$ large enough
\begin{align*}
  |\tau(0, v) - \mu(v)| & \leq  |\tau(0, v) - \tau(0, w(\|v\|_2))| + |\tau(0, w(\|v\|_2)) - \tau(0, \|v\|_2 u)|  \\
  &\quad \ + |\tau(0, \|v\|_2 u) - \mu(\|v\|_2 u)| + |\mu(\|v\|_2 u) - \mu(v)|.
\end{align*}
We bound each summand separately. By definition of $\tau'$ and $w(r)$ and the fact that $v$ and $w(r)$ are contained in a spherical shell of thickness $\rho$, we have that
\begin{align*}
  |\tau(0, v) - \tau(0, w(\|v\|_2))| & \leq \tau( v, w(\|v\|_2))
  \leq \tau'( v, w(\|v\|_2)) \\
  & \leq \lambda \|v - w(\|v\|_2)\|_2
  \leq \lambda \|v\|_2 c_4(\rho, \delta)
\end{align*}
where $c_4(\rho, \delta)$ is a function that tends to zero, when $\rho$ and $\delta$ tend to zero (see Lemma \ref{le:app2}). The same argument yields that
\[
  |\tau(0, w(\|v\|_2)) - \tau(0, \|v\|_2 u)| \leq \lambda \|v\|_2 c_4(\rho, \delta).
\]
Due to \eqref{eq:lipschitzbound_on_seminorm}
\begin{align*}
  |\mu(\|v\|_2 u) - \mu(v)| &\leq \max_{i \in [d]} \mu(\mathbf{e}_i) \| \|v\|_2 u - v\|_2 \leq \max_{i \in [d]} \mu(\mathbf{e}_i) \|v\|_2 c_4(\rho, \delta).
\end{align*}
Combining these bounds, we obtain that for every $\varepsilon > 0$ and direction $u \in \cS^{d-1}$ there are parameters $\lambda, r_0 \in \R$ and $\rho, \delta > 0$ such that a.s.\ for all $v \in S(u, \infty , \delta) \cap \Z^d$ with $\|v\|_2 > r_0$
\[
  \frac{1}{\|v\|_2} |\tau(0, v) - \mu(v)| < \varepsilon.
\]
Hence we have proven the convergence of \eqref{eq:shape_theorem_random_metric:uniform-limit} when $x$ is restricted to some cone intersected with $\Z^d$. The compactness of the sphere ensures, that for any given $\varepsilon > 0$ we find a finite number of cones with this convergence, that cover the whole space.\qed\bigskip

We already mentioned, that the problem with Theorem \ref{thm:shape_theorem_random_metric} is, that $\mu_\tau$ might be degenerate. In fact it is easy to show, that $\mu_\tau$ is a norm if and only if there is a shape $S_\tau$ such that a.s.\ for any $\varepsilon > 0$ and large enough $n \in \N$
\begin{align}\label{eq:shape-theorem-with-limit-shape-S}
  (1-\varepsilon) S_\tau \subset \frac{1}{n} B_n^\tau(Z_0) \subset (1 + \varepsilon) S_\tau.
\end{align}
In this case $S_\tau = \{x \in \R^d \mid \mu_\tau(x) \leq 1\}$ is a convex set with $S_\tau = -S_\tau$. We say \emph{$\tau$ satisfies a shape theorem with limit shape $S_\tau$}. It is already clear, that a weak shape theorem holds under rather weak conditions. However we want to find conditions that ensure a shape theorem with a limit shape.

%In it may sometimes be no ``real'' shape theorem. If $\mu_\tau$ is no norm but only a semi-norm then there are directions in which the space can be travelled at arbitrary high speed and hence the nice asymptotic behaviour
%\begin{align}\label{eq:real_shape_theorem}
%  (1-\varepsilon) S \subset \frac{B_t^\tau(0)}{t} \subset (1+ \varepsilon) S,
%\end{align}
%for $\varepsilon > 0$, $t$ large enough and $S := \{x \in \R^d \mid \mu_\tau(x) \leq 1\}$ fails.

The random pseudometrics that are induced by the ``\emph{fastest path mechanism}'' form a very large and important class of examples. We will inspect this class in more detail now.

Let $L$ be random time functional that maps continuous curves onto non-negative real numbers or infinity, i.e.\ $L(\gamma)$ tells us how long it takes to traverse the curve $\gamma$. This map shall be additive, i.e.\ $L(\gamma_1 \circ \gamma_2) = L(\gamma_1) + L(\gamma_2)$ for curves $\gamma_1, \gamma_2:[0,1] \to \R^d$ with $\gamma_1(1) = \gamma_2(0)$, where $\gamma_1 \circ \gamma_2$ is the concatenation of the two curves. Then the induced pseudometric $\tau_L$ is defined by
\begin{align}
  \tau_L(x,y) := \inf_{\gamma:x\leftrightarrow y} L(\gamma), \quad x,y \in \R^d
\end{align}
where the infimum is taken over all continuous curves from $x$ to $y$.

A very interesting example for a random pseudometric induced by such random time functional is the Riemannian FPP model from \cite{lagatta2009shape}. In a simple version of this model, a positive random field over $\R^d$ is generated and $L(\gamma)$ is just the integral of the random field along $\gamma$.

If we assumed $L(\gamma)$ to be larger than some fraction of the Euclidean length of $\gamma$ this would ensure that $\tau_L$ satisfies a shape theorem with limit shape, but it would be way too restrictive.

In fact, it should be possible to prove a shape theorem even if ``locally'' some curves are very short, as long as they behave nicely under a ``global'' perspective. The following lemma will make this precise. We define an auxiliary random field $W = (W_v)_{v \in \Z^d}$ depending on two parameters $\delta, \rho > 0$ via
\begin{align}
\begin{aligned}\label{eq:def-W}
  W_v &:= \ind\{\text{there is a $\gamma:x \leftrightarrow y$ with $x \in \delta (v + [-\tfrac{1}{2},\tfrac{1}{2}]^d)$ and}\\
  &\hspace{2cm} \text{$y \in \R^d \setminus \delta (v + [-\tfrac{3}{2}, \tfrac{3}{2}]^d)$ such that $L(\gamma) < \rho$} \}.
\end{aligned}
\end{align}
To shorten the notation we will write $v^{\square}$ for $\delta (v + [-\tfrac{1}{2},\tfrac{1}{2}]^d)$. Furthermore we have to introduce the set $\cA$ of connected subsets of $\cZ^d$ containing the origin. The elements of $\cA$ are called lattice animals and we denote by $\cA_n$ the animals consisting of $n \in \N$ vertices.

\begin{Lemma}\label{le:pos-time-const}
  Let $L$ be a random additive time functional such that $\tau_L$ is an ergodic pseudometric satisfying \eqref{eq:shape_theorem_random_metric:mom_cond}. If there are parameters $\delta, \rho > 0$ such that a.s.
  \begin{align}\label{eq:pos-time-const:animal-cond}
    c_1 := \limsup_{n \to \infty } \max_{\alpha \in \cA_n} \frac{1}{n} \sum_{v \in \alpha} W_v < 1,
  \end{align}
  then $\tau_L$ satisfies a shape theorem with limit shape $S_\tau$.
\end{Lemma}
\emph{Proof:} We only have to make sure, that $\mu(u) := \mu_{\tau_L}(u) > 0$ for all $u \in \cS^{d-1}$. For a given $u \in \cS^{d-1}$ and large enough $n$ there is a constant $c_2$ such that any curve from the origin to $nu$ intersects at least $c_2 \|nu\|_2$ different cubes $v^\square$ with $v\in \Z^d$. For a given curve $\gamma:0\leftrightarrow nu$ we denote the set of points $v \in \Z^d$ where $v^\square$ is intersected by $\alpha(\gamma)$ and observe that $\alpha(\gamma) \in \bigcup_{k \geq c_2\|nu\|_2} \cA_k$. Hence we find a subset $\alpha_0(\gamma) \subset \alpha(\gamma)$ of size at least $(1-c_1)|\alpha(\gamma)|$ of points $v$ such that $W_v = 0$. Moreover we find a subset $\alpha_1(\gamma) \subset \alpha_0(\gamma)$ such that $|\alpha_1(\gamma)| \geq |\alpha_0(\gamma)|/3^d$ and $\|v - w\|_\infty \geq 3$ for all $v \neq w \in \alpha_1(\gamma)$. The boxes $\delta (v + [-\tfrac{3}{2}, \tfrac{3}{2}]^d)$ and $\delta (w + [-\tfrac{3}{2}, \tfrac{3}{2}]^d)$ are disjoint for $v \neq w \in \alpha_1(\gamma)$. By the definition of $W$, the restriction of $\gamma$ to one of these boxes needs at least a time $\rho$ to be traversed. Combining this with the fact that $|\alpha_1(\gamma)| \geq c_2 (1 - c_1) 3^{-d} n$ shows, that the total time to traverse $\gamma$ is at least $\rho c_2 (1 - c_1) 3^{-d} n$ and hence $\mu(u) \geq c_2 (1 - c_1) 3^{-d} \rho > 0$.\qed \bigskip

The assumption of Lemma \ref{le:pos-time-const} is rather strict in general, but it holds in the case of Riemannian FPP (see \cite{lagatta2009shape}), if the random field inducing $L$ is $k$-dependent and satisfies some minor moment condition. The application of Theorem \ref{thm:shape_theorem_random_metric} and Lemma \ref{le:pos-time-const} is discussed in the following section.

\section{FPP on random tessellations}\label{sec:fpp-on-random-tess}
Theorem \ref{thm:shape_theorem_random_metric} and Lemma \ref{le:pos-time-const} may be applied to first passage percolation on random tessellations. Informally this model consists of a random tessellation, where a passage time is attached to each face of a cell. The time functional $L(\gamma)$ of a curve $\gamma$ just adds up the passage times of all faces that are intersected by $\gamma$. Formally we will use a marked particle process of faces with mark space $[0,\infty )$. This leads to the following definitions and notational conventions (a broad introduction to point processes and random tessellations can be found in \cite{schneider2008stochastic}).

Let $D$ be a metric space equipped with the Borel-$\sigma$-algebra $\cB(D)$. We write $\mathbf{N}(D)$ for the set of locally finite counting measures on $D$ and equip it with the $\sigma$-algebra $\cN(D)$ generated by the sets $\{\eta \in \mathbf{N}(D) \mid \eta(A) = k\}$, $A \in \cB(D), k \in \N_0$. A measure $\eta \in \mathbf{N}(D)$ is called locally finite iff $\eta(A) < \infty $ for all bounded $A \in \cB(D)$. A measurable mapping $\Phi:\Omega \to \mathbf{N}(D)$ is called a \emph{point process} on $D$ and is to be interpreted as a random collection of points in $D$. Each point process permits a representation
\[
  \Phi = \sum_{i = 1}^{\Phi(D)} \delta_{\zeta_i},
\]
where $\delta$ is the Dirac measure and $(\zeta_i)_{i \in \N}$ are $D$-valued random variables \cite[Lemma 3.1.3]{schneider2008stochastic}. The measure $\Theta := \E\Phi$ on $D$ is called the \emph{intensity measure} of $\Phi$. In the important special case, where $D = \R^d$ and $\Theta = \gamma \lambda^d$ we call $\gamma$ the \emph{intensity} of $\Phi$ ($\lambda^d$ is the Lebesgue measure on $\R^d$). Random tessellations will later be defined by letting $D$ be the space $\cC^d$ of compact and convex subsets of $\R^d$ equipped with the Hausdorff metric.

\begin{figure}
  \centering
  \begin{tikzpicture}
    \begin{scope}
      \clip(-4,-3.8) rectangle (4.5,3.8);
      \node at (0,0) {\includegraphics[scale=.6]{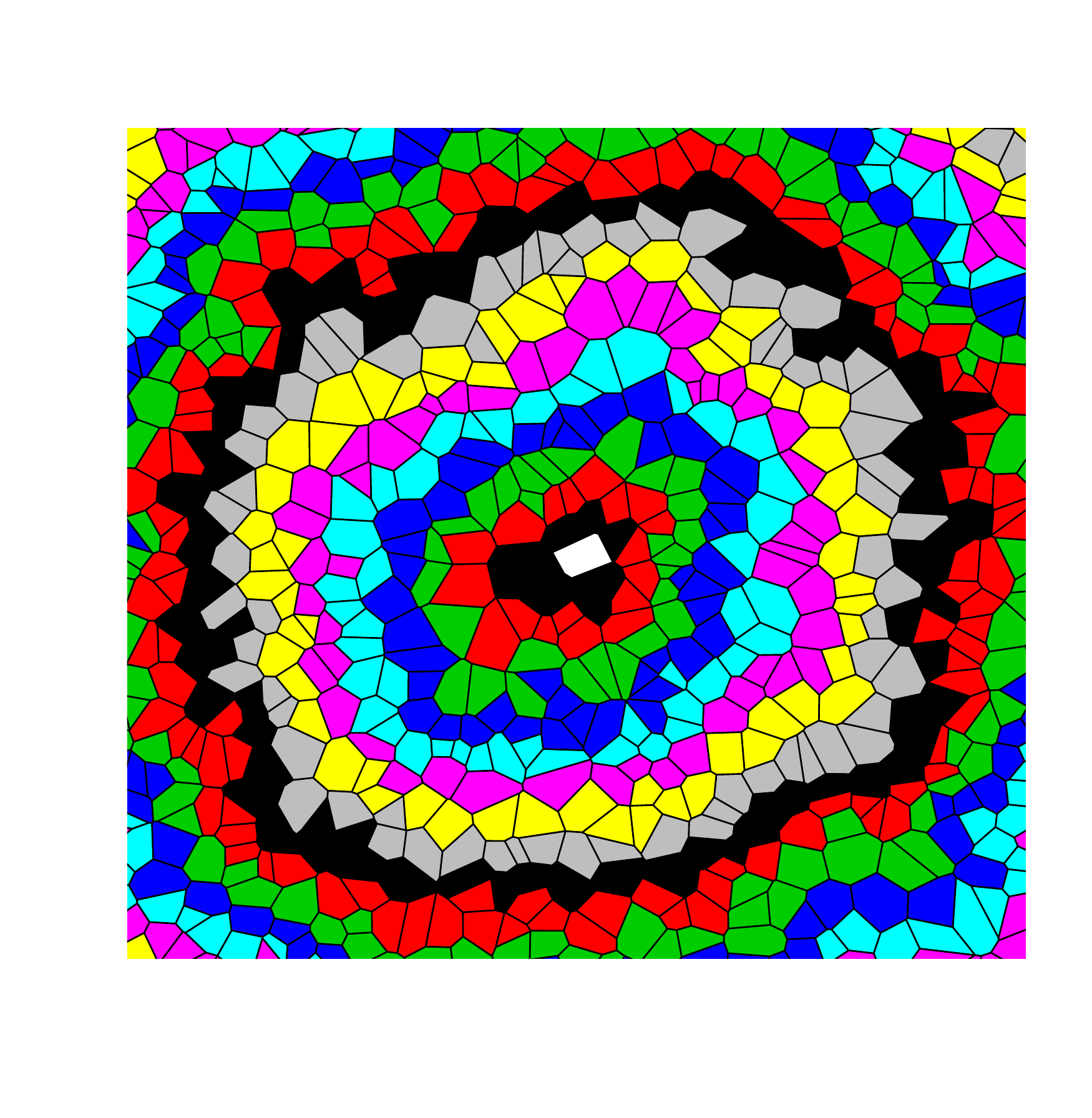}};
    \end{scope}
  \end{tikzpicture}
  \caption{A realization of a Voronoi tessellation $m$ induced by a determinantal point process. Points with the same distance to the origin with respect to $\tau_{m,\mathbf{1}}$ have the same color.}
  \label{fig:fpp-on-voronoi-of-dpp}
\end{figure}

To be able to add a passage time to each face of a cell, it is convenient to work with marked point processes. If we have a point process $\Phi = \{\zeta_1, \zeta_2, \dots\}$ on $D$ with representation $\sum_{i = 1}^{\Phi(D)} \delta_{\zeta_i}$ and a (possibly random) real valued sequence $X = (X_i)_{i \in \N}$ we call
\[
  \Phi_X := \sum_{i = 1}^{\Phi(D)} \delta_{(\zeta_i, X_i)}
\]
a \emph{marked version} of $\Phi$. In the case where $X$ is a $[0,\infty )$-valued i.i.d. sequence with marginal distribution $\Q$ that is independent of $\Phi$, we call $\Phi_X$ the \emph{independently marked version} of $\Phi$ with mark distribution $\Q$.

%An important special case is the one, where $X$ is an i.i.d. sequence with marginal distribution $\Q$ that is independent of $\Phi$. We call $\Phi_X$ the \emph{independently marked version} of $\Phi$ with mark distribution $\Q$ then.

From now on we will only work with point processes on $\R^d$ or $\cC^d$ and their marked versions (a marked point process on $D$ is also a point process on $D \times [0,\infty )$). If $D$ is equal to either of these spaces, the canonical translation operator $\theta_x:\mathbf{N}(D) \to \mathbf{N}(D)$, $x \in \R^d$ is defined by
\[
  \theta_x \eta(A) := \eta(A - x), \quad \eta \in \mathbf{N}(D), A \in \cB(D).
\]
We use the same notation for the shift $\theta_x:\mathbf{N}(D \times [0,\infty )) \to \mathbf{N}(D \times [0,\infty ))$ on the marked spaces defined by
\[
  \theta_x \eta(A \times B) := \eta( (A - x) \times B), \quad \eta \in \mathbf{N}(D), A \in \cB(D), B \in \cB([0,\infty )).
\]
This corresponds to the idea that only the points are shifted while each point retains its mark. A point process on $\cC^d$ is also called a \emph{particle process}. Stationarity and ergodicity are defined in the same way as in Section \ref{sec:shape-thm-ergodic-psmetrics}. It can be shown, that if $\Phi$ is ergodic, then its independently marked version is ergodic too; see \cite[Proposition 12.3.VI.]{daley2007introductionvolzwei}. We also remark, that in the important case, where the values of $X$ are a deterministic function of $\Phi$ that commutes with $\theta_x$ for all $x \in \R^d$, the marked version $\Phi_X$ is ergodic if $\Phi$ is ergodic.

%We recall, that a point process $\Phi$ is \emph{stationary} iff $\theta_x \Phi \overset{d}{=} \Phi$ for all $x \in \R^d$. Let $\cI$ be the $\sigma$-Algebra of translation invariant events, i.e.\ events $A \in \cN(D)$ with $T_x A = A$ for all $x \in \R^d$. A stationary point process $\Phi$ is called \emph{ergodic} if $\P[\Phi \in A] \in \{0,1\}$ for all $A \in \cI$.

%We recall the definition of the \emph{Laplace functional} of a point process $\Phi$ applied to a function $f:D \to [0,\infty )$
%\[
%  L_{\Phi}(f) := \E \exp\Big( -\int f(x)\ \Phi(dx)\Big).
%\]
%We will extend this definition to $\R$-valued functions $f$ and remark that the integral or the expectation might not exist in this case. An introduction to point processes and the Laplace functional for random measures can be found in \cite{daley2003introduction}.

After the introduction of point and particle processes, we turn to tessellations. A set $Z \in \cC^d$ with non-empty interior is called a \emph{cell}. A countable set $m := \{Z_1, Z_2, \dots\}$ of cells is called a \emph{tessellation} (or \emph{m}osaic) if
\begin{enumerate}
  \item each ball in $\R^d$ is intersected by at most a finite number of cells of $m$,
  \item the cells of $m$ cover $\R^d$,
  \item the interiors of any two distinct cells in $m$ doesn't overlap.
\end{enumerate}
The cell of $m$ that contains $x \in \R^d$ is denoted by $Z_x(m)$ (if there is more than one cell containing $x$ we chose an arbitrary rule to break ties). The cell $Z_0(m)$ is called the \emph{zero cell}.
%\begin{figure}
%  \centering
%    \includegraphics[scale=.9]{Bilder/G_m_and_G_m_star}
%    \caption{A section of a Poisson hyperplane tessellation $m$ and its induced graphs \textcolor{red}{$\cG_{\mathfrak{m}}$} und \textcolor{blue}{$\cG_{\mathfrak{m}}^*$}}
%  \label{fig:tessellation_and_induced_graphs}
%\end{figure}

Each tessellation $m$ is completely determined by the $d-1$-dimensional faces of its cells, i.e. by the set $\overline m := \{Z_1 \cap Z_2 \mid Z_1, Z_2 \in m, \text{dim}(\text{aff}(Z_1 \cap Z_2)) = d-1\}$ (aff$(\cdot)$ denotes the smallest affine space containing $\cdot$). Let $\mathbf{M} \subset \mathbf{N}(\cC^d)$ be the set of tessellations and let $\overline{\mathbf{M}} := \{\overline m \mid m \in \mathbf{M}\}$ be the set of face ensembles corresponding to a tessellation.

Observe, that any tessellation $m \in \mathbf{M}$ induces a graph $\cG_m := (m, E_m)$ with vertex set $m$. Two cells $Z_1, Z_2 \in m$ are adjacent in $\cG_m$ iff $Z_1 \cap Z_2 \in \overline m$ and hence we will identify $E_m$ with $\overline m$. The zero cell $Z_0(m)$ is the root $\mathbf{0}$ in $\cG_m$.

We denote by $\mathfrak{F}(\mathbf{M}) := \cN(\cC^d)|_\mathbf{M}$ the trace of $\cN(\cC^d)$ on $\mathbf{M}$ and call a measurable mapping $M:\Omega \to \mathbf{M}$ a random tessellation. Hence a random tessellation is a point process of convex compact particles that form a tessellation. In the same spirit let $\overline{\mathbf{M}}_t := \{(F_i, X_i)_{i \in \N} \mid (F_i)_{i \in \N} \in \overline{\mathbf{M}}, X_j \in [0,\infty ), j \in \N\}$ be the set of \emph{marked face ensembles} with the $\sigma$-Algebra $\mathfrak{F}(\overline{\mathbf{M}}_t) := \cN(\cC^d \times [0,\infty ))|_{\overline{\mathbf{M}}_t}$. For $\overline m = \{F_1, F_2, \dots\} \in \overline{\mathbf{M}}$ and a random or deterministic $[0, \infty )$-valued sequence $X = \{X_1, X_2, \dots\}$ we define the marked face ensemble $\overline m_X := \{(F_1, X_1), (F_2, X_2), \dots\} \in \overline{\mathbf{M}}_t$. %A marked tessellation $m_X$ induces a colored graph $\cG_{m,X} := (\cG_{m}, c)$ where $c(Z_i) := X_i$, $i \in \N$ (the same notations are used for a random tessellations $M$ in place of $m$).

Each marked face ensemble $\overline m_X = \{(F_1, X_1), (F_2, X_2), \dots\}$ induces the time functional $L_{m,X}$ via
\begin{align}
  L_{m,X}(\gamma) := \sum_{i \in \N} \ind\{F_i \cap (\gamma \setminus \gamma(1)) \neq \emptyset \} X_i,
\end{align}
where $\gamma:[0,1] \to \R^d$ is a curve. It is only for technical reasons that we exclude $\gamma(1)$ in the definition (e.g. to have additivity of $L$ for concatenated curves). The corresponding pseudometric $\tau_{L_{m,X}}$ is denoted by $\tau_{m,X}$.

At first we want to find a simple condition on $m$ and the mark distribution, such that we may apply Theorem \ref{thm:shape_theorem_random_metric} to $\tau_{m,X}$.

\begin{Lemma}\label{le:shape-theorem-rand-tess}
  Let $\overline M_X$ be an independently marked version of the face ensemble $\overline M$ of an ergodic random tessellation $M$. If there is an $\varepsilon > 0$ such that
  \begin{align}\label{eq:shape-theorem-rand-tess:mom-cond1}
    \E \Bigg[ \bigg( \sum_{Z \in M} \ind\{Z \cap [-1,1]^d \neq \emptyset \} \bigg)^{d+\varepsilon} \Bigg] < \infty
  \end{align}
  and
  \begin{align}\label{eq:shape-theorem-rand-tess:mom-cond2}
    \E[X_1^{d+\varepsilon}] < \infty,
  \end{align}
  then $\E[ \max \left\{ \tau_{M,X}(0,x) \mid \|x\|_\infty \leq 1 \right\}^{d+\varepsilon}] < \infty$ and $\tau_{M,X}$ satisfies a weak shape theorem.
\end{Lemma}
\emph{Proof:} The random face ensemble $\overline M$ is ergodic as $\theta_x \overline{M} = \overline{\theta_x M}$ and hence $\overline M_X$ is ergodic too as it is an independently marked version of $\overline M$. The pseudometric $\tau_{M,X}$ is a deterministic function of $\overline M_X$ that commutes with $\theta_x$ for each $x \in \R^d$, i.e.\ $\theta_x \tau_{M,X} = \tau_{\theta_x M, X}$. This yields the ergodicity of $\tau_{M,X}$.
Let $N := \{Z \in M \mid Z \cap [-1,1]^d \neq \emptyset \}$ be the set of cells that intersect the cube $[-1,1]^d$. The induced subgraph of $\cG_{M}$ with vertex set $N$ is a.s.\ connected as a.s.\ the intersection of every cell with the unit cube is either $d$-dimensional or empty. Hence there is a spanning tree of this graph with edge set $K$ of size $|N| - 1$. This implies that
\[
  \max\{\tau_{M,X}(0,x) \mid \|x\|_{\infty } \leq 1\} \leq \sum_{i \in \N} \ind\{F_i \in K\} X_i \overset{d}{=} \sum_{i = 1}^{|N|-1} X_i.
\]
Applying Jensen's inequality to this sum yields
\[
  \E \Bigg[  \bigg( \sum_{i = 1}^{|N|} X_i \bigg)^{d+\varepsilon}  \Bigg] \leq \E \bigg[ |N|^{d+\varepsilon-1} \sum_{i = 1}^{|N|} X_i^{d+\varepsilon} \bigg] = \E[|N|^{d+\varepsilon}] \E[X_1^{d+\varepsilon}]
\]
and hence the moment condition on $\tau_{M,X}$.\qed\bigskip

We remark, that the proof of \ref{le:shape-theorem-rand-tess} may be adapted to much more complicated time functionals. On could for example think of a model that assigns a speed to each cell and the time to traverse a curve could be defined as the Euclidean length weighted with the inverse speed of the corresponding cell plus the time needed to pass through the cell boundaries. Also deterministic passage times are possible. For example the passage time of a cell face could be a function of its $d-1$ Hausdorff measure. All these constructions would preserve the ergodicity. The moment conditions would have to be adjusted accordingly, of course.

Our concern now is, to tell when $\tau_{M,X}$ satisfies a shape theorem with limit shape $S_\tau$. To check condition \eqref{eq:pos-time-const:animal-cond} we will recall the notion of a \emph{tame tessellation} introduced in \cite{ziesche2016bernoulli}. We fix a grid width $\delta > 0$ like in the introduction of $W$ (see \eqref{eq:def-W}) and recall the notation $\zeta^{\square } := \delta ( \zeta + [-\tfrac 12, \tfrac 12]^d )$ for $\zeta \subset \Z^d$.

For a stationary random tessellation $M$ two auxiliary random fields $Y := (Y_v)_{v \in \Z^d} := (Y_v(M, \delta))_{v \in \Z^d}$ and $U := (U_v)_{v \in \Z^d} := (U_v(M, \delta))_{v \in \Z^d}$ were defined by
\begin{align}
  \begin{aligned}
    Y_v & := |\{Z \in M \mid z(Z) \in v^{\square \delta}\}| \\
    U_v & := \ind\{\text{a cell of $M$ intersects $v^{\square \delta}$ and $\{w \in \Z^d \mid \|w - v\|_\infty \geq 2\}^{\square \delta}$}\}.
  \end{aligned}
\end{align}

\begin{Def}
  An ergodic random tessellation $M$ is called \emph{tame} if there is a $\delta > 0$ such that
  \begin{enumerate}
    \item[(T1)] there is a $c_1 \in \R$ such that
        \[
          \limsup_{n\to \infty } \max_{\alpha \in \cA_n^{(d)}} \frac{1}{n} \sum_{v \in \alpha} Y_v \leq c_1,
        \]
    \item[(T2)] there is a $c_2 < 1$ such that
        \[
          \limsup_{n\to \infty } \max_{\alpha \in \cA_n^{(d)}} \frac{1}{n} \sum_{v \in \alpha} U_v \leq c_2.
        \]
  \end{enumerate}
\end{Def}

The most basic, but nevertheless interesting question might be, what happens in the case of a trivial distribution of passage times, i.e.\ when $X = \mathbf{1} := (1,1,\dots)$ a.s.\ . This corresponds to the usual graph metric on $\cG_M$. If there is a curve $\gamma$ connecting $x \in \R^d$ with $y \in \R^d$ such that $L_{M,\mathbf{1}}(\gamma) < 1$, then $x$ and $y$ have to be contained in the same cell. Hence $U_v \geq W_v$ for $v \in \Z^d$ and therefor condition (T2) implies \eqref{eq:pos-time-const:animal-cond}. This ensures that $\tau_{M,\mathbf{1}}$ satisfies a shape theorem with limit shape if the moment condition \eqref{eq:shape-theorem-rand-tess:mom-cond1} holds.

In the case of i.i.d.\ random passage times $X$, we would like to find a feature of the mark distribution, that characterizes the existence of a limit shape. This is at moment not in reach, as we don't know enough about the percolation properties of random tessellations in this general setup. However we can show, that a limit shape exists, if (T1) holds and the mark distribution has no large atom at zero.

\begin{Lemma}
  Let $\overline M_X$ be an independently marked version of the face ensemble $\overline M$ of an ergodic random tessellation $M$ such that the moment conditions \eqref{eq:shape-theorem-rand-tess:mom-cond1} and \eqref{eq:shape-theorem-rand-tess:mom-cond2} hold. If $M$ is a tame random tessellation then there is a constant $c_1 > 0$ such that $\tau_{M,X}$ satisfies a shape theorem with limit shape if $\P[X_1 = 0] < c_1$.
\end{Lemma}
\emph{Proof:} According to Theorem 5.2. from \cite{ziesche2016bernoulli} the tameness of $M$ implies, that there is a constant $c_2$ such that a.s.\ $|\cA_n(\cG_M)| \leq c_2^n$ for large enough $n$ where $\cA_n(\cG_M)$ is the set of connected subgraphs of $\cG_M$ with $n$ vertices one of which being the zero cell $Z_0(M)$ (this is in analogy to $\cA_n$ for $\cZ^d$). Let $m \in \mathbf{M}$ be such that $|\cA_n(\cG_m)| \leq c_2^n$ for all $n$ larger than some $n_0$.

For some constant $c_3 > 0$ and $n \in \N$ we consider the event $E_n$ that holds if there is a curve $\gamma$ starting in the origin and traversing exactly $n$ cell boundaries of $m$ such that $L_{m,X}(\gamma) \leq c_3 n$. The curve $\gamma$ corresponds to a path $\tilde \gamma$ in $\cG_m$ that traverses the same faces $F_i \in m$. We identify $\tilde \gamma$ with the sequence of indices of faces traversed by $\gamma$ and apply Markov's inequality to obtain that
\begin{align*}
  \P[E_n] & \leq \sum_{\tilde \gamma} \P\bigg[  \sum_{i \in \tilde \gamma} X_i \leq c_3 n \bigg] \\
  & \leq c_2^{n+1} e^{t c_3 n} \E[e^{-tX_1}]^{n} \\
  & \leq c_2^{n+1} e^{t c_3 n} (\P[X_1 < \varepsilon] + e^{-t \varepsilon} \P[X_1 \geq \varepsilon])^{n}
\end{align*}
for all $t > 0$, where the second inequality is due to the fact that each path starting in $Z_0(m)$ is also an animal of $\cG_m$.

If $\P[X_1 = 0] < c_2^{-1}$ we find $\varepsilon, t, c_3 > 0$ (in this order), such that $e^{t c_3} (\P[X_1 < \varepsilon] + e^{-t \varepsilon} \P[X_1 \geq \varepsilon]) < c_2^{-1}$. Hence by the Borel-Cantelli lemma a.s.\ only a finite number of the events $E_n$ hold for this $c_3$. This implies that a.s.\ $L(\gamma) \geq c_3 n$ for all $n$ large enough and any curve $\gamma$ that traverses at least $n$ cell boundaries.

Due to the considerations before this lemma, we already know, that there is a constant $c_4 > 0$ such that a.s.\ for large enough $n$ any curve intersecting $n$ cell boundaries is contained in the Euclidean ball of radius $c_4 n$. This finishes the proof as for $x \in \R^d$ eventually $L(\gamma) \geq c_3/c_4 \|x\|_2$ if $\gamma:0 \leftrightarrow x$ and $\|x\|_2$ tends to infinity.\qed\bigskip

In \cite{ziesche2016bernoulli} there are various examples of tame tessellations most notably Voronoi tessellations induced by determinantal point processes or Poisson cluster processes.

Summarizing this section we showed three things.
\begin{itemize}
  \item Under rather weak conditions, a weak shape theorem holds for $\tau_{M,X}$.
  \item If in addition (T2) holds, then $\tau_{M,\mathbf{1}}$ satisfies a shape theorem with limit shape.
  \item If $M$ is even tame, i.e. (T1) and (T2) hold, then $\tau_{M,X}$ satisfies a shape theorem with limit shape, if $\P[X_1 = 0] < c_1$ for some $c_1 > 0$.
\end{itemize}

\section{A spatial ergodic theorem on random graphs}\label{sec:spatial-erg-thm-on-graphs}
Our knowledge about FPP on a random tessellation $M$ yields insights on the spatial structure of the graph $\cG_M$. In particular we want to show in this section, that the multivariate ergodic theorem of Wiener \cite[Theorem 10.14]{kallenberg2002foundations} that considers averages over a growing sequence of fixed convex observation windows, maybe modified such that the average is taken over \emph{random balls} in $\cG_M$.

A random tessellation $M$ is called \emph{weakly ergodic} if $\P[M \in A] \in \{0,1\}$ for all $A \in \bigcap_{x \in \R^d} \cI_x$. This notion is called ergodic in the literature most of the time but the differences to what we call ergodic are rather small \cite{pugh1971ergodic}. We recall Wiener's theorem in the case where it is applied to a weakly ergodic random tessellation.

\begin{NamedThm}\label{thm:spatial-ergodic-thm}(Wiener)\\
  Let $M$ be a weakly ergodic random tessellation and $W_1 \subset W_2 \subset \dots$ a sequence of bounded convex subsets of $\R^d$ with inner radii tending to infinity. Then for any measurable $f:\mathbf{M} \to [0,\infty )$
  \begin{align}
    \lim_{n\to \infty }\frac{1}{\Vol(W_n)} \int_{W_n} f(\theta_x M) \ dx = \E[f(M)] \quad \text{a.s.\ .}
  \end{align}
\end{NamedThm}

Our goal is to apply this Theorem in the following informal setting. Let $h$ be a translation invariant function of a cell and the surrounding tessellation. Then we want to show, that for an ergodic tessellation $M$ the limit
\begin{align}\label{eq:informal-average-over-ball}
  \lim_{n\to \infty } \frac{1}{|B_n^M|} \sum_{Z \in B_n^M} h(Z,M)
\end{align}
a.s.\ exists and is equal to the expected value of $h(Z,M)$ under the Palm probability measure of $M$. The set $B_n^M := B_n^{\cG_M}(Z_0(M))$ in the equation is the ball of radius $n$ around the zero cell in the graph metric on $\cG_M$. To this end we have to introduce some basics of Palm calculus. A broader introduction may be found in \cite{schneider2008stochastic}.

A function $z:\cC^d \to \R^d$ is called a \emph{center function} if it is measurable and translation covariant, i.e.\ if $z(Z + x) = z(Z) + x$ for $x \in \R^d$, $Z \in \cC^d$. For the rest of the paper $z$ will be an arbitrary center function such that $z(Z) \in Z$. One might think of $z$ for example as the barycenter of $Z$ (further examples can be found in \cite{schneider2008stochastic}).

When dealing with Palm probabilities, we will often have to access the centers of cells in our random tessellation. Therefor it is useful to interpret a random tessellation $M = \{Z_i \mid i \in \N\}$ as a marked point process $M' := \{(z(Z_i), Z_i - z(Z_i)) \mid i \in \N\}$ of cell centers marked with the corresponding cell.

We define the \emph{cell intensity}
\begin{align}
  \gamma_M := \E[|\{(x,Z) \in M' \mid x \in [0,1]^d\}| ] = \E[|\{Z \in M \mid z(Z) \in [0,1]^d\}|]
\end{align}
and the \emph{Palm distribution}
\begin{align}
  \P_{M'}^0[A] := \frac{1}{\gamma_M} \E\bigg[ \sum_{(x,Z) \in M'} \ind\{x \in [0,1]^d\} \ind\{(Z, \theta_x M')\}\bigg],
\end{align}
where $A \in \cB(\cC^d) \otimes \cN(\R^d \times \cC^d)$. The distribution $\P_{M'}^0$ may be understood as the joint distribution of the zero cell and $M$ conditioned on the event that one of the points of $M'$ lies in the origin. The corresponding cell is then called the \emph{typical cell}. Informally we can think of the typical cell as a cell that was picked ``uniformly'' from all cells of $M$. Theorem \ref{thm:ergodic-theorem-on-graphs} will not only establish an ergodic theorem on graphs but will also clarify how Palm probabilities can be interpreted.

For a measurable function $f:\R^d \times \cC^d \times \mathbf{N}(\R^d \times \cC^d) \to [0,\infty )$ we recall the well known Campbell Theorem (see \cite[Theorem 3.5.3]{schneider2008stochastic})
\begin{align}
  \E\bigg[ \sum_{(x,Z) \in M'} f(x,Z,M') \bigg] = \gamma_M \iint f(x,Z,\theta_x m)\ \P_{M'}^0(d(Z,m))\ dx.
\end{align}

For the rest of this section, we want to study a translation invariant function $h:\cC^d \times \mathbf{N}(\R^d \times \cC^d) \to [0,\infty )$, i.e.\ $h(Z,m') = h(Z + x, \theta_{-x} m')$ for all $x \in \R^d$. By applying Campbell's Theorem to the function $f$ defined by $f(x,Z,m') := \ind\{-x \in Z\} h(Z_0(m), m') / \Vol(Z)$ we further improve our understanding of Palm probabilities as
\begin{align}
\begin{aligned}\label{eq:palm-expectations}
  \E\left[ \frac{h(Z_0(M), M')}{\Vol(Z_0(M))}\right] & = \E\bigg[ \sum_{(x,Z) \in M'} f(x,Z,M') \bigg] \\
  & = \gamma_M \iint h(Z_0(m), m') \frac{\ind\{-x \in Z\}}{\Vol(Z_0(m))}\ dx\ \P_{M'}^0(d(Z,m')) \\
  & = \gamma_M \int h(Z, m') \ \P_{M'}^0(d(Z,m')) \\
  & = \gamma_M \E_{M'}^0[h(Z_0(M), M')],
\end{aligned}
\end{align}
where $\E_{M'}^0$ is the expectation with respect to $\P_{M'}^0$. Equation \eqref{eq:palm-expectations} implies that the expected value of a functional of the zero cell and its surrounding tessellation is equal to the cell intensity times the expected value of this functional weighted with the volume of the zero cell under the Palm probability. In particular
\begin{align}\label{eq:cell-intensity}
  \E[\Vol(Z_0(M))^{-1}] = \gamma_M.
\end{align}

%We say a random tessellation $M$ \emph{satisfies a shape theorem with limit shape $S$} if for any $\varepsilon > 0$ and large enough $n$
%\begin{align}\label{eq:shape-theorem-with-limit-shape-S}
%  (1-\varepsilon) S \subset \frac{1}{n} B_n^\tau(Z_0) \subset (1 + \varepsilon) S,
%\end{align}
%where $\tau := \tau_{M,\mathbf{1}}$.
If $\tau$ satisfies a shape theorem with limit shape, then the set of points that can be reached by crossing at most $n$ cell boundaries looks approximately like $n S_\tau$. Hence the idea is, to use $n S_\tau$ as $W_n$ in Wieners theorem and then replace them by $B_n^\tau$.

In the following theorem we write $B_n^M$ for the ball $B_n^{\cG_M}(Z_0(M)) \subset M$ in the graph $\cG_M$, $B_r^\tau$ for $B_r^{\tau_{M,\mathbf{1}}}(0) \subset \R^d$ and $B_r$ for $B_r(0)$. We want to point out that $B_n^\tau$ is equal to $\bigcup_{Z \in B_n^M} Z$.

\begin{Thm}\label{thm:ergodic-theorem-on-graphs}
  Let $M$ be a weakly ergodic random tessellation such that $\tau_{M, \mathbf{1}}$ satisfies a shape theorem with limit shape $S := S_{\tau_{M, \mathbf{1}}}$. Then
  \begin{align}\label{eq:ergodic-theorem-on-graphs:claim1}
    \lim_{n\to \infty } \frac{|B_n^{M}|}{n^d} = \gamma_M \Vol(S) \quad a.s.
  \end{align}
  and
  \begin{align}\label{eq:ergodic-theorem-on-graphs:claim1.1}
    \lim_{n\to \infty } \frac{|B_n^{\tau}|}{n^d} = \Vol(S) \quad a.s.\ .
  \end{align}
  Furthermore for any measurable $g:\mathbf{M} \to [0,\infty )$
  \begin{align}\label{eq:ergodic-theorem-on-graphs:claim2}
    \lim_{n\to \infty } \frac{1}{\Vol(B_n^\tau)} \int_{B_n^\tau} g(\theta_{x} M) \ dx = \E[g(M)] \quad a.s.\ .
  \end{align}
  Finally for a translation invariant $h:\cC^d \times \mathbf{N}(\R^d \times \cC^d) \to [0,\infty )$
  \begin{align}\label{eq:ergodic-theorem-on-graphs:claim3}
    \lim_{n\to \infty } \frac{1}{|B_n^M|} \sum_{Z \in B_n^M} h(Z,M) = \E_{M'}^0[h(Z_0,M')] \quad a.s.\ .
  \end{align}
\end{Thm}
\emph{Proof:} We start with \eqref{eq:ergodic-theorem-on-graphs:claim2}. For $\varepsilon \in (0,1)$ applying \eqref{eq:shape-theorem-with-limit-shape-S} and Theorem \ref{thm:spatial-ergodic-thm} with $f := g$ and $W_n := n(1- \varepsilon) S$ yields that
a.s.
\begin{align*}
  \E[f(M)] & = \lim_{n\to \infty } \frac{1}{\Vol(n (1 - \varepsilon) S)} \int_{n(1-\varepsilon) S} f(\theta_x M) \ dx \\
  & \leq \lim_{n\to \infty } \frac{(1+\varepsilon)^d (1-\varepsilon)^{-d}}{\Vol(B_n^\tau)} \int_{B_n^\tau} g(\theta_x M) \ dx.
\end{align*}
The lower bound works in the same way and the third assertion follows.

To prove \eqref{eq:ergodic-theorem-on-graphs:claim1} and \eqref{eq:ergodic-theorem-on-graphs:claim1.1} we observe, that $Z_0(\theta_x M) = Z_{x}(M) - x$ and hence
\begin{align}\label{eq:ergodic-theorem-on-graphs:2}
  \int_{B_n^\tau} \Vol(Z_0(\theta_x M))^{-1} \ dx & = \sum_{Z \in B_n^M} \int_{Z} \Vol(Z_x(M))^{-1} \ dx = |B_n^M|.
\end{align}
The arguments leading to the third assertion also suffice to show, that
\begin{align}\label{eq:ergodic-theorem-on-graphs:3}
  \lim_{n\to \infty } \frac{1}{\Vol(nS)} \int_{B_n^\tau} g(\theta_{x} M) \ dx = \E[g(M)] \quad a.s.\ .
\end{align}
Applying this to the function $g(M) := \Vol(Z_0(M))^{-1}$ as well as using \eqref{eq:ergodic-theorem-on-graphs:2} and \eqref{eq:cell-intensity} yields, that
\[
  \lim_{n\to \infty } \frac{|B_n^M|}{n^d \Vol(S)} = \E[\Vol(Z_0(M))^{-1}] = \gamma_M.
\]
Using \eqref{eq:ergodic-theorem-on-graphs:claim2} instead of \eqref{eq:ergodic-theorem-on-graphs:3} in this argument yields \eqref{eq:ergodic-theorem-on-graphs:claim1.1}.

Finally, we notice that due to the translation invariance of $h$
\begin{align*}
  \sum_{Z \in B_n^M} h(Z,M) & = \sum_{Z \in B_n^M} \int_{Z} \frac{h(Z_x(M), M)}{\Vol(Z_x(M))} \ dx \\
  & = \sum_{Z \in B_n^M} \int_{Z} \frac{h(Z_0(\theta_{x} M), \theta_{x} M)}{\Vol(Z_0(\theta_{x} M))} \ dx = \int_{B_n^\tau} g(\theta_x M) \ dx,
\end{align*}
where $g(m) := h(Z_0(m), m)/\Vol(Z_0(m))$. Combining \eqref{eq:ergodic-theorem-on-graphs:claim2} applied to this $g$ with \eqref{eq:ergodic-theorem-on-graphs:claim1}, \eqref{eq:ergodic-theorem-on-graphs:claim1.1} and \eqref{eq:palm-expectations} finishes the proof.\qed \bigskip

A first non-trivial application of Theorem \ref{thm:ergodic-theorem-on-graphs} is obtained, when we define $h(Z,M)$ as the probability that $Z$ lies in the boundary of two infinite clusters in the Bernoulli percolation model on $M$ as defined in \cite{ziesche2016bernoulli}. The existence of the point-wise limit was an assumption in Theorem 2.2 in \cite{ziesche2016bernoulli} that showed the uniqueness of the infinite cluster in the Bernoulli percolation model on $M$.

The simple proof of Theorem \ref{thm:spatial-ergodic-thm} makes clear, that in the same way one could average over random balls induced by other random pseudometrics that satisfy a shape theorem with limit shape. However our focus was to improve understanding of global properties of random tessellations and \eqref{eq:ergodic-theorem-on-graphs:claim3} helps a lot with this.

\section{Poisson hyperplane tessellations}\label{sec:poisson-hyperplane-tess}
\begin{figure}
  \centering
  \begin{tikzpicture}
    \begin{scope}
      \clip(-3,-2.5) rectangle (3.5,2.5);
      \node at (0,0) {\includegraphics[scale=.31]{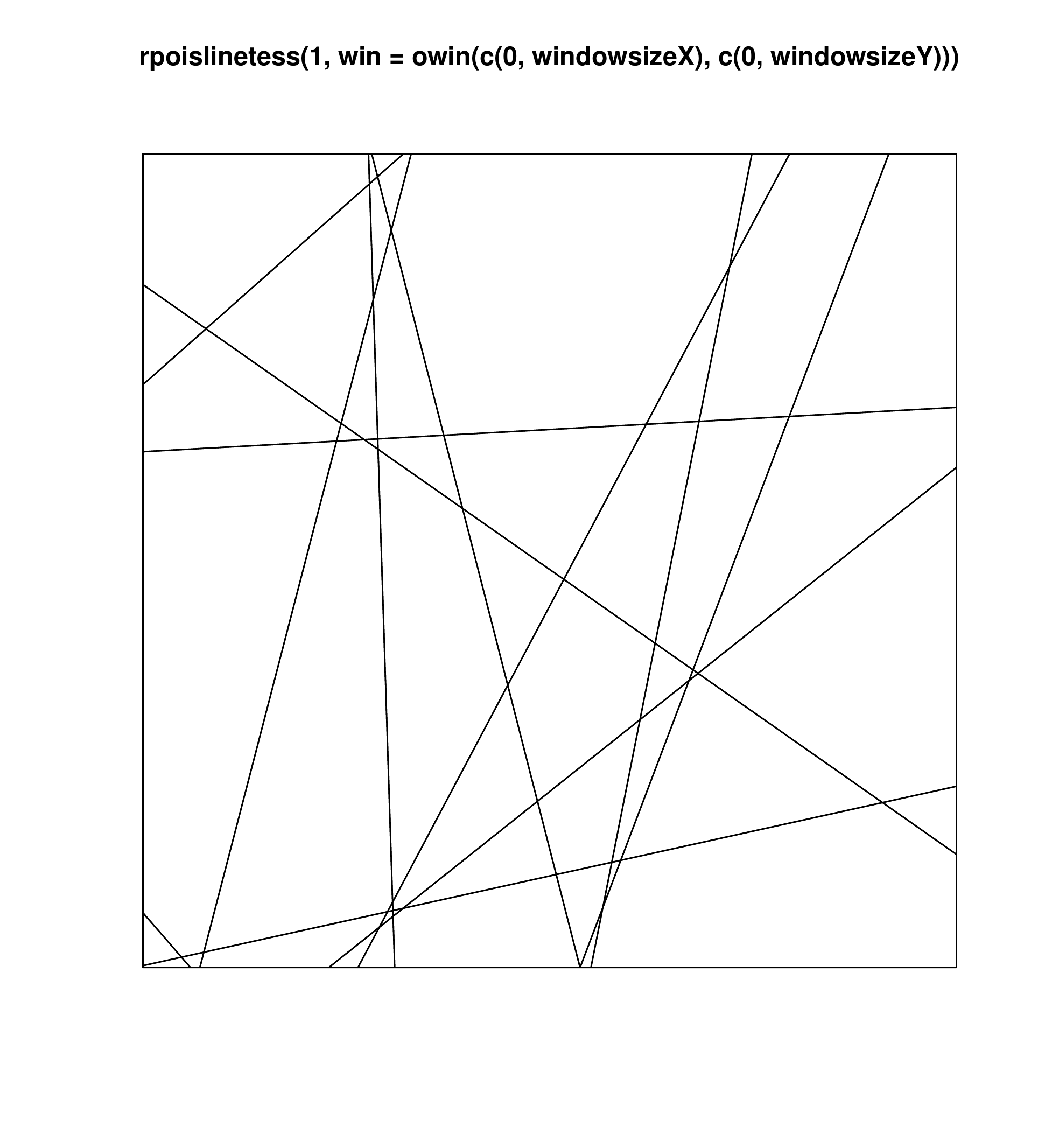}};
    \end{scope}
  \end{tikzpicture}
  \caption{A realization of a Poisson hyperplane tessellation where $\varphi$ is the uniform distribution.}
  \label{fig:poisson-line-tess}
\end{figure}

The Poisson hyperplane tessellation (PHT) is generated by an infinite set of hyperplanes, that is distributed randomly in $\R^d$; see Figure \ref{fig:poisson-line-tess}. While the Voronoi tessellation is tame for certain classes of point processes \cite{ziesche2016bernoulli}, it was argued, that this is not the case for the PHT. However in contrast to the Bernoulli percolation model, the FPP is so well suited to the PHT, that not only are we able to show that the corresponding pseudometric satisfies a shape theorem, but we are in fact able to \emph{determine the limit shape} in certain cases. Moreover we give bounds on the speed of convergence for the limit in \eqref{eq:shape_theorem_random_metric:uniform-limit}.

We will introduce the PHT formally first. For a direction $u \in \cS^{d-1}$ and an $r \geq 0$ let $E(u,r) := \{x \in \R^d \mid \langle x, u\rangle = r\}$ be the hyperplane perpendicular to $u$ with a distance $r$ to the origin.

We consider a Poisson process $\Phi$ on $\cS^{d-1} \times [0,\infty )$ with intensity measure $\E \Phi = \gamma \varphi \otimes \lambda^d$, where $\gamma \in (0,\infty )$, $\varphi$ is an even probability measure on $\cS^{d-1}$ that is not concentrated on a great subsphere and $\lambda^d$ is the Lebesgue measure. This is the point process characterized by the two properties that first $\Phi(A)$ is Poisson distributed with parameter $\gamma \varphi \otimes \lambda^d(A)$ for any measurable $A \subset \cS^{d-1} \times [0,\infty )$ and second $\{\Phi(A_i) \mid i \in I\}$ is a family of independent random variables for pairwise disjoint $A_i$ with $i$ in some finite index set $I$.

We identify $\Phi$ with the set $\{E(u, r) \mid (u,r) \in \Phi\}$ of random hyperplanes and call $\Phi$ a Poisson hyperplane process. The assumptions on $\varphi$ imply that $\Phi$ is ergodic (even mixing \cite[Theorem 9.3.6]{schneider2008stochastic}) and a.s.\ tessellates the space into bounded cells. In particular $\Phi$ induces a random tessellation $M := M(\Phi)$ with a law that is determined by $\gamma$ and $\varphi$. We call $M$ the \emph{Poisson hyperplane tessellation with directional distribution $\varphi$ and intensity $\gamma$}. A nice introduction to PHT can be found in \cite{schneider2008stochastic}.

A central object to describe a PHT is the so-called \emph{associated zonoid}. That is a convex set with support function
\begin{align}
  h(x) := \gamma \int_{\cS^{d-1}} |\langle x, u\rangle |\ \varphi(du), \quad x \in \R^d.
\end{align}
At the same time
\begin{align}
  h(x) = \frac{1}{2} \E[ |[0,x] \cap \Phi| ], \quad x \in \R^d,
\end{align}
where $[0,x]$ is the line from $0$ to $x$ \cite[Section 4.4]{schneider2008stochastic}.

Now there are three natural ways to assign passage times to the cell boundaries of $M$. The first way is to assign a constant time of $1$ to each boundary. The second way would be to assign i.i.d. random passage times to the hyperplanes of $\Phi$ that would be inherited by each cell boundary in that hyperplane. The last way would be to take i.i.d. passage times for all cell faces of $\overline M$.

For the third option we obtain the conditions for Lemma \ref{le:shape-theorem-rand-tess} after observing, that the number of cells intersecting the unit cube is smaller than two to the number of hyperplanes intersecting it. This number is Poisson distributed and hence the moment condition \eqref{eq:shape-theorem-rand-tess:mom-cond1} is satisfied, yielding a weak shape theorem for $\overline M_X$ if $X_1$ has a finite $d + \varepsilon$ moment. It is an open question, under which conditions on the law of $X_1$, there is a limit shape in this case.

The second way is considerably easier to handle, such that we are able to give the exact distribution of $\tau(0,x)$ in this case and compute the limit shape $S_\tau$. For the rest of this Section, let $\Phi = \{E_1, E_2,\dots\}$ be a Poisson hyperplane process and $X = \{X_1, X_2,\dots\}$ the usual i.i.d. sequence of passage times. We call $\Phi_X := \{(E_1, X_1), (E_2, X_2), \dots\}$ the timed Hyperplane process, define the time functional $L_{\Phi,X}$ via
\[
  L_{\Phi,X}(\gamma) := \sum_{i \in \N} \ind\{(\gamma \setminus \gamma(1)) \cap E_i \neq \emptyset \} X_i,
\]
where $\gamma:[0,1] \to \R^d$ is a curve and obtain the corresponding pseudometric $\tau := \tau_{L_{\Phi,X}}$.

The crucial observation for this model is, that \emph{the direct way is always the fastest}. If two points $x,y \in \R^d$ lie on different sides of a hyperplane $E$, each curve connecting $x$ and $y$ has to cross $E$. The straight line $[x,y]$ crosses only such hyperplanes and hence is as fast as possible. This implies that
\begin{align}\label{eq:pht-tau-is-compound-poisson}
  \tau(x,y) = \sum_{i \in \N} \ind\{E_i \cap [x,y) \neq \emptyset \} X_i,
\end{align}
which is a compound Poisson distribution with an expected value of
\[
  \E[\tau(0,x)] = \E[X_1] \E[|[0,x] \cap \Phi|] = 2h(x) \E[X_1].
\]
We may also immediately conclude that
\[
  \max\{t(0,x) \mid \|x\|_\infty  \leq 1\} \leq \sum_{i \in \N} \ind\{E_i \cap [-1, 1]^d \neq \emptyset \} X_i
\]
has a finite $d+\varepsilon$ moment for some $\varepsilon > 0$ as long as the $d+\varepsilon$ moment of $X_1$ is finite, which we will always assume in this Section.

Hence we know, that $\tau$ satisfies at least a weak shape theorem. However the equations \eqref{eq:shape_theorem_random_metric:1} and \eqref{eq:pht-tau-is-compound-poisson} show, that
\[
  \mu(x) = \inf_{n \in \N} \frac{\E[\tau(0,nx)]}{n} = 2 \E[X_1] h(x).
\]
This quantity is larger than zero for all $x \neq 0$ if $\E[X_1] > 0$ as $\varphi$ is not degenerate, which implies that $\tau$ satisfies a shape theorem with limit shape $S = \{x \in \R^d \mid 2 \E[X_1] h(x) \leq 1\}$.

Knowing the exact limit shape and the distribution of $\tau(0,x)$ produces the question, whether we can say something about the speed of convergence in the limit considered in Theorem \ref{thm:shape_theorem_random_metric}. We will answer this question with the following theorem, where we restrict ourselves to the case $X = \mathbf{1}$. The result holds also for random $X$, but the speed will  depend a lot on the tail behavior of $X_1$. It is however very easy to incorporate this into the proof by replacing \eqref{eq:pht-speed-of-conv:11} and \eqref{eq:pht-speed-of-conv:12} with a concentration result for the corresponding compound Poisson variable.

\begin{Thm}\label{thm:pht-speed-of-conv}
  Let $X = \mathbf{1}$ and $\Phi$ be a Poisson hyperplane process with directional distribution $\varphi$, intensity $\gamma$ and corresponding norm $\mu := \mu_\tau$. There is a constant $c_1 \in \R$ such that for all $\varepsilon > 0$ and $r > 0$
  \begin{align}\label{eq:pht-speed-of-conv:claim}
    \P[\exists\ x \in r\cS^{d-1} \colon |\tau(0,x) - \mu(x)| > \varepsilon r] \leq c_1 \varepsilon^{-2(d-1)} \exp\left( -\frac{r \varepsilon^2}{m}\right),
  \end{align}
  where $m := 8 \max_{u \in \cS^{d-1}} \mu(u)$.
\end{Thm}
\emph{Proof:} We work as in the proof of Theorem \ref{thm:shape_theorem_random_metric} and quantify the different details. On the one hand, we have to consider how many cones are needed to cover the whole space (this geometrical argument was done in the appendix Lemma \ref{le:app3}) on the other hand, we need a concentration inequality for the difference between $\tau$ and $\mu$ in a spherical section of a given size.

Let $\varepsilon > 0$. For a given direction $u \in \cS^{d-1}$, radius $r \in [0,\infty ]$ and opening parameter $\delta \in (0,1]$ we recall the definition of a spherical section
\begin{align}
  S(u,r,\delta) := \{x \in \R^d \mid \|x\|_2 \leq r,\ \langle \tfrac{x}{\|x\|_2}, u\rangle \geq 1 - \delta\}.
\end{align}
We define the set
\begin{align}
  \cE_1(u, r, \delta) := \{E(v,a) \mid \exists\ \tilde u \in \cS^{d-1},\ \langle u, \tilde u\rangle \geq 1 - \delta \colon \langle \tilde u, v\rangle \geq \tfrac{a}{r}\}
\end{align}
of hyperplanes that intersect $S(u,r,\delta)$ and the set
\begin{align}
  \cE_2(u, r, \delta) := \{E(v,a) \mid \forall\ \tilde u \in \cS^{d-1},\ \langle u, \tilde u\rangle \geq 1 - \delta \colon \langle \tilde u, v\rangle \geq \tfrac{a}{r}\}
\end{align}
of hyperplanes that intersect each straight line in $S(u,r, \delta)$ from $0$ to $x$ with $\|x\|_2 = r$. This definition implies, that
\begin{align}\label{eq:pht-speed-of-conv:4}
  \max\{\tau(0,x) \mid x \in S(u,\infty , \delta),\ \|x\|_2 = r\} \leq \Phi(\cE_1(u, r, \delta))
\end{align}
and
\begin{align}
  \min\{\tau(0,x) \mid x \in S(u,\infty , \delta),\ \|x\|_2 = r\} \geq \Phi(\cE_2(u, r, \delta)).
\end{align}

For fixed $u \in \cS^{d-1}$ and $r \in (0,\infty )$, the expected number of hyperplanes in $\cE_i(u,r,\delta)$ should converge to $\mu(ru)$ as $\delta \to 0$. We will calculate how fast this happens. For $i \in \{1,2\}$ it is easily seen, that $E(v,a) \in \cE_i(u,r,\delta)$ iff $E(v,\frac{a}{r}) \in \cE_i(u,1, \delta)$ and hence
\begin{align}
  \E[\Phi(\cE_i(u,r, \delta))] = \gamma \int \int_0^\infty  \ind_{\cE_i(u, r, \delta)}(E(v,a))\ da\ \varphi(dv) = r \E[\Phi(\cE_i(u,1,\delta))].
\end{align}
Furthermore we have, that for any $u,v,w \in \cS^{d-1}$ and $\delta \in (0,1]$ with $\langle u,w \rangle \geq 1 - \delta$
\begin{align}
  \langle v, w \rangle \leq \langle u, v \rangle + \delta + \sqrt{8\delta} \leq \langle u,v \rangle + c_2 \sqrt{\delta},
\end{align}
where $c_2 := 1 + \sqrt{8}$ (see Lemma \ref{le:app1}). This yields, that
\begin{align}
\begin{aligned}\label{eq:pht-speed-of-conv:8}
  &\quad \ \E[\Phi(\cE_1(u, 1, \delta))] \\
  & = \gamma \int \int_{0}^\infty \ind\{\exists\ w \in \cS^{d-1},\ \langle w, u\rangle \geq 1- \delta \colon \langle w, v \rangle \geq a\}\ da\ \varphi(dv)\\
  & \leq \gamma \int \int_{0}^\infty \ind\{\langle u,v\rangle + c_2 \sqrt{\delta} \geq a \}\ da\ \varphi(dv)\\
  & \leq \gamma \int \int_{0}^\infty \ind\{\langle u,v\rangle \geq a \}\ da\ \varphi(dv) + \gamma c_2 \sqrt{\delta} \\
  & = \E[\Phi(\{E(v,a) \mid E(v,a) \cap [0,u] \neq \emptyset \})] + \gamma c_2 \sqrt{\delta} \\
  & = \mu(u) + \gamma c_2 \sqrt{\delta}
\end{aligned}
\end{align}
and
\begin{align}
\begin{aligned}
  &\quad \ \E[\Phi(\cE_2(u, 1, \delta))] \\
  & = \gamma \int \int_{0}^\infty \ind\{\forall\ w \in \cS^{d-1},\ \langle w, u\rangle \geq 1- \delta:\ \langle w, v \rangle \geq a\}\ da\ \varphi(dv)\\
  & \geq \gamma \int \int_{0}^\infty \ind\{\langle u,v\rangle - c_2 \sqrt{ \delta} \geq a \}\ da\ \varphi(dv)\\
  & \geq \mu(u) - \gamma c_2 \sqrt{\delta}.
\end{aligned}
\end{align}

In the next step we will quantify how much $\mu$ might vary within $r \cS^{d-1} \cap S(u,\infty ,\delta)$. Due to \eqref{eq:lipschitzbound_on_seminorm} and Lemma \eqref{le:app2}
\begin{align}\label{eq:pht-speed-of-conv:10}
  |\mu(x) - \mu(ru)| \leq \max_{i \in [d]} \mu(\mathbf{e}_i) \|x - ru\|_1 \leq c_3 r \sqrt{\delta}
\end{align}
for $x \in S(u,\infty ,\delta)$, $\|x\|_2 = r$ and $c_3 := 2\sqrt{2} \max_{i \in [d]} \mu(\mathbf{e}_i)$.

Before putting everything together, we recall the well known bounds
\begin{align}\label{eq:pht-speed-of-conv:11}
  \P[P \geq \E[P] + x] \leq \exp\left( \frac{-x^2}{2\E[P]}\right), \quad x \geq 0
\end{align}
and
\begin{align}\label{eq:pht-speed-of-conv:12}
  \P[P \leq \E[P] - x] \leq \exp\left( \frac{-x^2}{2\E[P]}\right), \quad x \geq 0
\end{align}
for a Poisson distributed random variable $P$.

Let $\delta \in (0,1]$. Due to Lemma \ref{le:app3} there is a $c_4 \in \R$, a $k \leq c_4 \delta^{1-d}$ and a set of directions $\{u_i \in \cS^{d-1} \mid i \in [k]\}$ such that $\cS^{d-1}$ is covered by the spherical segments $S(u_i, 1, \delta)$ and hence $\R^d$ is covered by the one-sided cones $S(u_i, \infty ,\delta)$. Combining this with  \eqref{eq:pht-speed-of-conv:10} and \eqref{eq:pht-speed-of-conv:4} yields
\begin{align*}
  & \quad \ \P[\exists\ x \in r \cS^{d-1} \colon \tau(0,x) - \mu(x) > \varepsilon r] \\
  & \leq \sum_{i \in [k]} \P[\max\{\tau(0,x) \mid x \in S(u_i,\infty ,\delta),\ \|x\|_2 = r\} > r(\mu(u_i) - c_3 \sqrt{\delta} + \varepsilon)] \\
  & \leq \sum_{i \in [k]} \P[\Phi(\cE_1(u_i, r, \delta)) > r(\mu(u_i) - c_3 \sqrt{\delta} + \varepsilon)].
\end{align*}
Applying \eqref{eq:pht-speed-of-conv:11}, \eqref{eq:pht-speed-of-conv:8} and the fact that $k \leq c_4 \delta^{1-d}$, we obtain
\begin{align*}
  \P[\exists\ x \in r \cS^{d-1} \colon \tau(0,x) - \mu(x) > \varepsilon r] \leq \frac{c_4}{\delta^{d-1}} \exp\left( - \frac{r(\varepsilon - (\gamma c_2 + c_3)\sqrt{\delta})^2}{2 \mu(u_i)}\right).
\end{align*}
Choosing $\sqrt{\delta} = \varepsilon/2(\gamma c_2 + c_3)$ yields one half of \eqref{eq:pht-speed-of-conv:claim}. The other half is obtained in the very same way using all the lower bounds derived in the proof.\qed \bigskip

\section*{Acknowledgements}
This article covers parts of the results of the authors PhD thesis. The author wants to thank G\"unter Last for his support during this time and for the multitude of fruitful discussions.
%To show (T2) one might use \cite{liggett1997domination}

\section{Appendix}
We need some estimates about scalarproducts of spherical vectors.

\begin{Lemma}\label{le:app1}
  Let $\cH$ be a Hilbert space with inner product $\langle \cdot, \cdot\rangle$, then for $x,y,z \in \cH$
  \begin{align*}
    \langle x,y\rangle & \geq \langle x,z\rangle + \langle y,z\rangle - \langle z,z\rangle \\
    & \quad \quad - \sqrt{(\langle x,x\rangle + \langle z,z\rangle - 2 \langle x,z\rangle)(\langle y,y\rangle + \langle z,z\rangle - 2 \langle y,z\rangle)}.
  \end{align*}
  In particular, if $\langle x,x\rangle = \langle y,y\rangle = \langle z,z\rangle = 1$, then
  \[
    \langle x,y\rangle \geq \langle x,z\rangle + \langle y,z\rangle - 1 - 2 \sqrt{(1-\langle x,z\rangle) (1 - \langle y,z\rangle)}.
  \]
\end{Lemma}
\emph{Proof:} We have a norm on $\cH$ defined by $\|x\| := \sqrt{\langle x,x\rangle}$ and a metric on $\cH$ defined by $d(x,y) := \|x-y\| = \sqrt{\langle x,y\rangle + \langle y,y\rangle - 2 \langle x,y\rangle}$. Squaring both sides of the triangle inequality for this metric yields the assertion.\qed \bigskip

We recall the definition of a spherical sector $S(u,r, \delta) := \{x \in \R^d \mid \|x\|_2 \leq r, \langle \frac{x}{\|x\|_2}, u \rangle \geq 1 - \delta \}$, $u \in \cS^{d-1}, r \in [0,\infty ], \delta \in [0,1]$.

\begin{Lemma}\label{le:app2}
  Let $u \in \cS^{d-1}$, $0 \leq r_1 < r_2 < \infty $ and $\delta > 0$. For $x,y \in S(u,r_2, \delta) \setminus S(u, r_1, \delta)$
  \[
    \|x - y\|_2 \leq r_2 - r_1 + 2 r_2 \sqrt{2\delta}.
  \]
\end{Lemma}
\emph{Proof:} Let $z \in \{\lambda x \mid \lambda \geq 0\}$ such that $\|z\|_2 = \|y\|_2$. Then $z \in S(u,r_2, \delta) \setminus S(u, r_1, \delta)$ and $\|x - z\|_2 \leq r_2 - r_1$. Furthermore
\[
  \|z -\|z\|_2 u\|_2 = \sqrt{\langle z -\|z\|_2 u,z -\|z\|_2 u\rangle} \leq \sqrt{2\delta} \|z\|_2 \leq \sqrt{2 \delta} r_2
\]
and the same holds for $\|y - \|y\|_2 u\|_2$. Applying the triangle inequality finishes the proof.\qed\bigskip

\begin{Lemma}\label{le:app3}
  There is a constant $c_1 \in \R$ that depends only on the dimension $d$ such that for any $\delta > 0$ there is a set of directions $\{u_i \in \cS^{d-1} \mid i \in [k]\}$ such that
  \[
    \cS^{d-1} \subset \bigcup_{i \in [k]} S(u_i,1,\delta)
  \]
  and
  \[
    k \leq c_1 \delta^{1-d}.
  \]
\end{Lemma}
\emph{Proof:} We will construct the covering for fixed $1 \geq \delta > 0$. To this end we cover $\cS^{d-1}$ with cubes $Q_v := \frac{\delta}{2\sqrt{d}} ( [0,1)^d + v)$, $v \in \Z^d$ and name the set of cubes $U$. The cubes in $U$ have a diameter of $\delta /2$, hence no cube of $U$ can intersect the sphere $(1- \delta) \cS^{d-1}$ or $(1 + \delta) \cS^{d-1}$. This implies, that the union of cubes has a volume less or equal to $\kappa_d( (1 + \delta)^d - (1 - \delta)^d)$ where $\kappa_d$ is the volume of the unit ball in $\R^d$. This bounds the size of $U$ by
\[
  |U| \leq \frac{\kappa_d ( (1 + \delta)^d - (1 - \delta)^d)}{(2\sqrt{d})^{-d} \delta^d} \leq \frac{c_1}{\delta^{d-1}},
\]
with $c_1$ chosen appropriately.

For each $Q_v \in U$ we fix a $u_v \in Q_v \cap \cS^{d-1}$. Let $x \in Q_v$, then there is an $y$ such that $x = y + u_v$ and $\|y\|_2 \leq \delta/2$. This implies, that $\|x\|_2 \leq 1 + \delta/2$ and
\[
  \langle u, x \rangle = 1 + \langle u, y\rangle \geq 1 - \delta/2.
\]
Hence $\langle x/\|x\|_2, u\rangle \geq 1 - \delta$ and $x \in S(u_v, 1, \delta)$. We conclude that the spherical sectors $S(u_v, 1, \delta)$ with with $v$ such that $Q_v \in U$ yields a covering of $\cS^{d-1}$ that has the claimed size.\qed \bigskip

\bibliography{C:/Arbeit/LaTeX/myBib}
\end{document}